\documentclass[12pt]{amsart}
\usepackage{amscd,amssymb,righttag}
\newtheorem{thm}{Theorem}[section]
\newtheorem{corol}[thm]{Corollary}
\newtheorem{lemma}[thm]{Lemma}
\newtheorem{prop}[thm]{Proposition}
\theoremstyle{definition}
\newtheorem{defin}[thm]{Definition}

\theoremstyle{remark}
\newtheorem{remark}[thm]{Remark}
\newtheorem{example}[thm]{Example}
\numberwithin{equation}{section}
\font\smc=cmcsc10 at 12pt
\font\smallsmc=cmcsc10
\let\o=\operatorname
\let\cal=\mathcal
\def\gcoor#1#2,#3#4{#1^1,\dots,#1^#2,#3^1,\dots,#3^#4}

\def\iso{\kern.35em{\raise3pt\hbox{$\sim$}\kern-1.1em\to}\kern.3em}
\def\oskip{\par\vbox to4mm{}\par}

\def\D{{\cal D}}
\def\X{{\frak X}}\def\E{{\cal E}_{\infty}}
\def\Y{{\frak Y}}
\def\Cal #1{{\mathcal #1}} 
\def\R{{\Bbb R}}
\def\Z{{\Bbb Z}}

\def\N{{\Bbb N}}
\def\bysame{$\raise.2em\hbox to 3em{\hrulefill}$\thinspace, }
\begin{document}
\rightline{\today}
\oskip\oskip
\begin{center}
{\bf WHAT IS SUPERTOPOLOGY? $^\ddag$}
\par\oskip\oskip
{\smc Ugo Bruzzo}\thinspace \S,\ddag \ \ and \
{\smc Vladimir Pestov\thinspace \P}
\par\bigskip
{\S\thinspace Scuola Internazionale Superiore di Studi Avanzati (SISSA),}
\par {Via Beirut 2-4, 34014 Trieste, {\smc Italy}}
\par\smallskip
{\ddag\thinspace Dipartimento di Matematica, Universit\`a degli Studi di
Ge\-no\-va,}
\par{Via Dodecaneso 35, 16146 Ge\-no\-va, {\smc Italy}}
\par\smallskip
{\P\thinspace School of Mathematical and Computing Sciences,
Victoria University}
{\par of Wellington, P.O. Box 600, Wellington, {\smc New Zealand}}
\par\medskip{E-mail:
{\tt bruzzo@sissa.it},
{\tt vova@mcs.vuw.ac.nz}}
\par\medskip
{Home (V.P.): {\tt http://www.vuw.ac.nz/$^\sim$vova}}
\end{center}
\par\oskip
\begin{quote}\footnotesize
{\smallsmc Abstract.\ } We discuss the problem of finding an analogue of
the concept of a topological space in supergeometry, motivated by
a search for a procedure to compactify a supermanifold along
odd coordinates. In particular, we examine the
topologies
naturally arising on the sets of points of locally ringed superspaces,
and show that in the presence of a nontrivial odd sector such topologies
are never compact.
The main outcome of our
discussion is that not only the usual framework of supergeometry (the theory of
locally ringed spaces), but the more general approach of the
functor of points, need to be further enlarged.
\end{quote}
\renewcommand{\thefootnote}{\fnsymbol{footnote}}
\footnotetext[3]{Prepublished as Research Report RP-98-25, School of
Mathematical and Computing Sciences, Victoria University of Wellington,
October 1998.}
\renewcommand{\thefootnote}{\arabic{footnote}}
\oskip\section{Introduction}
Geometries with anticommuting variables (supergeometries) have  been introduced
in connection with several issues in theoretical physics, notably to study
supersymmetric field theories; physical motivations to introduce such
geometries are briefly discussed in the next section. Supergeometries have been
quite intensively studied in the 70s and 80s. (It is impossible to provide here
any exhaustive bibliography; we only cite \cite{BBH1} as a general
reference, and the basic works
\cite{Be2,BL1,Ko,L1,DW,R1,JP,Rt1,BBHP2}. More detailed
references will be provided later on.)

What are the proper analogues of major
concepts of topology for such geometries?
In particular, how can we find `topologies' which
are capable of carrying information about the structure of superspaces `in
the odd
directions'? This is not an idle question, as,
for example, finding the right `super' analogue of compactness
and the ways to compactify supermanifolds are likely to have an impact on the
formulation of some physical theories.
And the lack of a satisfactory cohomology
theory for
superspaces is just another manifestation of our failure to conceive
(super)geometric
objects which exhibit nontrivial topological structure in their odd sector.

This article is especially written for a topological audience
and aims at inviting researchers in topology to join the quest
for methods to `superise' their area of knowledge.

We begin with an outline of the origins and basic ideas of supergeometry
and its place in our days' theoretical physics.
Then we present the sheaf-theoretic setting of supermanifold theory,
serving as the basis for supergeometry, whose key notion is that
of a locally ringed superspace. Our presentation is essentially
self-contained.
We discuss the problem, first explicitely stated by
Leites and Manin, of finding supergeometric
analogues of compactness.

In supergeometry, to every superspace there is associated a 
covariant functor
from the category, $\cal G$, of all finite dimensional Grassmann
algebras and graded-preserving algebra homomorphisms, to the
category ${\cal S}ets$ of all sets and mappings. 
Such functors, $\X$, are termed {\it virtual superspaces.} 
The image of a Grassmann algebra, $\wedge(q)$, of rank $q$ under
a virtual superspace functor, $\X$, is denoted by
$pt_q(\X)$ and called the set of $q$-points of $\X$.
Let  ${\cal T}op$ denote the category of all topological spaces and
continuous mappings.
It is quite natural to take as a basic concept of
`supertopology' that of a {\it virtual topological superspace,} that is,
an object of the category
 ${\cal T}op^{\cal G}$, formed by all covariant functors from
$\cal G$ to  ${\cal T}op$ and the corresponding functorial
morphisms. We will
show that if $\frak X$ is a locally ringed superspace and the ground
field $k$ is topological (as is the case in all the standard examples),
then for every $q\in\N$ the set of $q$-points of $\frak X$ carries
a natural topology.
Therefore, a virtual superspace associated to an arbitrary
locally ringed superspace has in fact a richer structure -- that of
a virtual topological superspace.
We will show, however,
that a virtual topological superspace determined by a locally ringed
superspace and having a nontrivial odd sector is
never `supercompact.' At best, such spaces are locally compact.

The category of virtual topological superspaces possesses natural
compactifications;
however, it seems that some of the most interesting conjectural objects of
`supertopology,' such as the hypothetical `purely odd projective space,'
do not correspond to objects of this category.

In Section \ref{g-spaces}
we establish an isomorphism between the category of
all virtual topological superspaces and a certain full subcategory of
the category of topological $G$-spaces
for some special semigroup $G=\E$, thus (at least, formally) reducing the
theory of virtual topological superspaces to abstract topological dynamics and
also achieving a somewhat greater generality than that offered by the
functor of points.
Whether or not the dynamical approach offers any
tangible advantages, or at least provides a new vantage point from which
to survey the development of the theory and decide upon further
directions, remains to be seen.

We also briefly comment upon the status of cohomology theories for superspaces.

Our resuls therefore suggest that, in order to embrace the phenomenon of
compactness, the existing framework has to be further extended.
The answer to the question asked in the title of this article is
thus `nobody seems to know (yet)!'

\section{Supergeometry and physics}

The origins of supergeometry lie in theoretical physics, and are to be
sought for in the procedure of `integrating over fermion variables' in quantum
field theory \cite{Ma}. This operation was performed by a formal device, now
called `Berezin integration' \cite{Be1},
which has been given a precise meaning
in supermanifold theory in the papers \cite{HM2} and \cite{Rt2}.
The usage of anticommuting variables has been advocated also in
connection with the dynamics of classical spinning particles \cite{BM},
the theory of superintegrable systems \cite{Ku},
and  the BRS analysis of quantum field theories \cite{BPT1,BPT2,HQRD}.
However, the most relevant motivation for supergeometry is nowadays provided by
supersymmetric field theory and superstring theory.  For same samples of
applications of supergeometry to supersymmetric gauge theory and supergravity
cf.~\cite{BC,Br1,BL3}.

A supersymmetric field theory is a quantum field theory
involving both bosonic and fermionic  fields which is invariant under a
transformation which mixes the two types of fields. The first examples of such
theories are due to Volkov and Akulov
\cite{VA} and Wess and Zumino \cite{WZ} (for an introduction to
supersymmetry the reader may refer to \cite{W}). We offer now a brief
introduction to  classical (pre-quantized) supersymmetric field theory.

Let us start by recalling a few basic notions in   field theory.
The arena where physical facts take place is {\it spacetime,} which
is usually supposed to be four-dim\-en\-sion\-al, three of its dimensions
accounting
for the usual three spatial dimensions, and the fourth dimension being
identified with time. (However, in most supersymmetric and string
theories spacetime is assumed to be higher dimensional, usually
of dimension 10 or 11.)
 One usually regards spacetime (at least when no
gravitational forces are taken into account) as a Euclidean four-dimensional
differentiable manifold equipped with a flat pseudo-Riemannian metric of
signature (3,1) (so in  pseu\-do\--car\-tesian
coordinates $\{x^i\}$ the metric $g$ has the standard form
$g=\o{diag}(1,1,1,-1)$).\footnote{In theoretical physics the signature (1,3)
is more often used.}

There are
two types of elementary particles, bosons and fermions. The basic microscopic
constituents of matter are fermions, such as electrons, quarks,
neutrinos; they are represented in terms of fields by spinor fields on the
spacetime manifold $M$. In mathematical terms, this means that fermionic fields
are sections of   spin bundles associated to the (principal) bundle of
orthonormal
frames on $M$. On the other hand, bosons are the carriers (quanta) of the
forces
acting between elementary particles: so photons are responsible for
electromagnetic interactions, the $W$ and $Z$ particles for weak nuclear
interactions, and the gluons for strong nuclear interactions. Bosons are
represented by tensor fields (basically because they are mathematically
described
by connections on principal bundles whose base manifold is spacetime $M$).

At the classical (non-quantum) level all fields are supposed to satisfy partial
differential equations,  called {\it field equations,} which be expressed as
Euler-Lagrange equations associated with a suitable action functional.
The latter is usually assumed to be  {\it local},
i.e., to be expressible as an integral of the type
$$\int_M\,L\,d\mu\,,$$
where $L$ (the {\it Lagrangian function\/}) is a function of the fields and
their first derivatives with respect to the spacetime coordinates, and $d\mu$
is a suitable measure.

A hint to the usage of supergeometry in supersymmetric field theory may be
provided  by the Wess-Zumino model.  The fields  in this model are two
complex scalar fields $A$,
$F$, and  a Dirac spinor field, $\psi^\alpha$, $\alpha=1\dots 4$ (so $A$,
$F$ are bosonic fields, and $\psi$ is a fermionic field).
The Lagrangian of the model is (letting $\partial_i=\frac\partial{\partial
x^i}$)
$$L=-\tfrac12\sum_{i=1}^4\left(\partial_iA\,\partial^iA^\ast+
i\bar\psi\gamma^i
\partial_i\psi\right)+\tfrac12 FF^\ast,$$
where $^\ast$ denotes complex conjugation and a bar denotes taking
the adjoint spinor. Moreover,
$\partial^i=\sum_{k=1}^4g^{ik}\,\partial_k$, where $g$ with upper indices is
the inverse matrix to $g$, and the $\gamma^i$ (the Dirac matrices) are some
matrices which are related to the construction of the spin bundles \cite{Law}.

The supersymmetry transformations rules for this model are given by
\begin{align*} A&\mapsto A+i\bar\varepsilon\psi \\
\psi&\mapsto\psi+\partial_iA\,\gamma^i\varepsilon+F\varepsilon \\
F&\mapsto F+i\bar\varepsilon\gamma^i\partial_i\psi \,.
\end{align*}
These transformations   leave the Lagrangian $L$
invariant up to first order in $\varepsilon$ if the parameters
$\varepsilon^\alpha$ and the spinor components
$\psi^\alpha$ anticommute among
them,
$$\varepsilon^\alpha\varepsilon^\beta=-\varepsilon^\beta\varepsilon^\alpha,
\qquad
\varepsilon^\alpha\psi^\beta=-\psi^\beta
\varepsilon^\alpha,\qquad\psi^\alpha\psi^\beta=-\psi^\beta\psi^\alpha.$$

This simple model shows that already at the classical (non-quantum) level a
consistent formulation of a supersymmetric field theory requires  some
generalization of differential geometry which is able to
incorporate anticommuting objects. A first result in this
direction is due to Salam and Strathdee
\cite{SS}. Their construction is  purely local, and amounts
essentially to the description of the geometry of the tangent space to a
supermanifold at a fixed point, yet it contains many of the basic ideas. Salam
and Strathdee introduced the notion of {\it superspace,} formally
described as a space  parametrized by four real
coordinates
$\{x^i\}$ together with four additional coordinates $\{y^\alpha\}$ satisfying
the commutation rules
$x^iy^\alpha=y^\alpha x^i$,
$y^\beta y^\alpha=-y^\alpha y^\beta$.
A scalar field $\Phi(x,y)$ on superspace --- usually called a {\it superfield}
--- can be developed in powers of the `odd coordinates' $y^\alpha$,
\begin{align*}
\Phi(x,y)& =\,\phi_0(x)+\sum_{1\leq\alpha\leq 4}
y^\alpha\,\phi_\alpha(x)
+\sum_{1\leq\alpha<\beta\leq 4} y^\alpha y^\beta
\,\phi_{\alpha\beta}(x) \\
& +\sum_{1\leq\alpha<\beta<\gamma\leq 4} y^\alpha y^\beta y^\gamma
\,\phi_{\alpha\beta\gamma}(x)+y^1y^2y^3y^4
\,\phi_{1234}(x)\,.\end{align*}
The expansion is finite due to the nilpotency of the $y^\alpha$.
The coefficients of this expansion  can be expressed in terms of  the
fields of the Wess-Zumino model \cite{WB}.
The formulation of a supersymmetric field theory in terms of superfields
usually achieves  remarkable simplifications \cite{1001}.

A more intriguing  appearance of supermanifolds in theoretical physics is
the usage of {\it moduli spaces of super Riemann surfaces} in string theory.
The partition function of bosonic strings can be expressed in a form involving
integrals over the moduli spaces of Riemann surfaces of all genera \cite{GSW}.
There have been attempts to extend these results to the theory of superstrings
(strings with fermionic degrees of freedom), by devising a `superisation' of
Riemann surfaces, and studying the  moduli spaces of the resulting objects. Any
additional details on this topic would lead us too
astray; the interested reader may refer to \cite{Ma2,Fri,LR,FR}. One should
however notice that it is exactly in this context that a procedure to
compactify
supermanifolds also ``along the odd directions'' is likely to be more relevant.

\section{Graded algebra preliminaries}

Let $k$ denote an arbitrary field.
We will assume all algebras (over $k$) to be associative and
unital.
The word {\it graded} will be always synonymous with
${\Bbb Z}_2${\it -graded}.

A {\it graded vector space} is a vector space $E$ together with a
fixed direct sum decomposition $E=E_0\oplus E_1$. Elements of
$E_0$ are said to be {\it even,} while elements of $E_1$ are
referred to as {\it odd.} Correspondingly, the vector subspace
$E_0$ is called the {\it even part} ({\it sector}) of $E$, and
$E_1$ is the {\it odd part} ({\it sector}).
For any element $x\in E_0\cup E_1$ one denotes by
$\tilde x\in {\Bbb Z}_2$ the {\it parity} of $x$, determined by
the rule $x\in E_{\tilde x}$ and computed $\mod 2$. Notice that
elements of $E\setminus (E_0\cup E_1)$ have no parity. Elements of
$E_0\cup E_1$ are called {\it homogeneous.}

An example is supplied by the `arithmetical' graded
vector space $k^{m\vert n}$, which is just the vector space
$k^{m+n}$ equipped with the grading
$(k^{m\vert n})_0=k^m$, $(k^{m\vert n})_1=k^n$.

Most of the basic constructions of linear algebra extend to the
graded case and, in particular, direct sums and tensor products
of graded vector spaces carry a natural
(and often self-evident) grading. For details, we
refer the reader to \cite{Ma3}.

A {\it graded algebra} is an algebra carrying a structure of a
graded vector space, $A=A_0\oplus A_1$. The two structures are required to
agree with each other in the sense that
for every $x,y\in A_0\cup A_1$ one has
\begin{equation}
\widetilde{xy}=\tilde x+\tilde y
\end{equation}
An associative unital graded algebra
$\Lambda=\Lambda_0\oplus\Lambda_1$ is called {\it graded commutative}
if for all $x,y\in\Lambda_0\cup\Lambda_1$ one has
\begin{equation}
xy=(-1)^{\tilde x\tilde y}yx
\end{equation}
This condition means that two even elements always commute with each other,
as well as an even and an odd element, while two odd elements
always anticommute.

In particular, every commutative unital associative algebra, $A$, equipped
with the trivial grading ($A_0=A$, $A_1=(0)$),
yields an example of a graded commutative algebra.

An ideal $I\subseteq\Lambda$ of a graded algebra $A$
 is called {\it graded} if $I=I_0\oplus I_1$,
where $I_i=I\cap\Lambda_i$, $i=1,2$.
A graded commutative algebra $\Lambda$ is called {\it local}
if it contains a unique maximal proper graded ideal $\frak m$; the quotient
$\Lambda/\frak m$ is a field, called the {\it residue field}
of $\Lambda$. It is always
a field extension of the ground field $k$, possibly a proper one.

The most important single example of a local graded commutative
algebra is provided by the exterior, or Grassmann, algebra $\wedge(q)$ of
rank $q$.
One way to describe $\wedge(q)$ is
as an associative unital
algebra freely generated by $q$ pairwise anticommuting elements
$\xi_1,\xi_2,\dots,\xi_q$. In other words, every element of
$\wedge(q)$ is a polynomial in the variables $\xi_i$ with coefficients
from $k$, having the form
\begin{equation}
\label{exp}
\sum_{\alpha\in 2^{\{1,\dots,q\}}} a_{\alpha}\xi^{\alpha},
\end{equation}
where $\alpha\subseteq \{1,2,\dots,q\}$,
$\xi^{\alpha}=\xi_{\alpha_1}\xi_{\alpha_2}
\dots\xi_{\alpha_{\vert\alpha\vert}}$, and $\xi^{\emptyset}=1$.
Generators are
subject to the anticommutation relations
\begin{equation}
\xi_i\xi_j=-\xi_j\xi_i\quad \mbox{for all $i,j=1,\dots,q$}
\end{equation}
Every collection of $q$ odd elements $\alpha_1,\dots,\alpha_q$ of
an arbitrary graded commutative (unital associative) algebra
$\Lambda$ determines a unique graded algebra homomorphism
from $\wedge(q)$ to $\Lambda$ with $\xi_i\mapsto\alpha_i$. Such a
homomorphism (or sometimes the corresponding collection 
$\alpha_1,\dots,\alpha_q\in\Lambda_1$
that fully determines it) is called, after the fashion of algebraic
geometry, a $\Lambda${\it -point} of $\wedge(q)$.

An element of $\wedge(q)$ is even (odd) if it can be represented
as a sum of monomials
in the free anticommuting generators $\xi_i$
having each  (respectively, odd) degree.
The maximal graded ideal of $\wedge(q)$ consists of all nilpotent
elements, which are exactly those polynomials $p(\xi)$
having vanishing constant term.

More generally, every Grassmann algebra $\wedge(q)$
supports a natural filtration
by a decreasing family of graded ideals $I_k$, $k=1,2,\dots,q$,
where $I_k$ is formed by all elements in whose expansion
(\ref{exp}) one has $a_{\alpha}=0$ whenever $\vert\alpha\vert<k$.
In particular, $I_1$ is the maximal graded ideal of $\wedge(q)$. It is easy to
verify that
$I_k\cdot I_m\subseteq I_{k+m}$.

The {\it augmentation homomorphism} (also called, in supergeometric
jargon, the
{\it body map}) is the quotient homomorphism
$\beta\colon\wedge(q)\to\wedge(q)/I_1\cong k$, associating to every
polynomial (\ref{exp}) the constant term $a_{\emptyset}$.

The simplest nontrivial example of a Grassmann algebra is
that of rank one,
$\wedge(1)$. As a vector space, it is the direct
sum of $k$ and the linear span of an odd generator $\xi$. Every
element of $\wedge(1)$ is then uniquely represented in the form
$a+b\xi$, where $a,b\in k$ and $\xi^2=0$. The structure of
$\wedge(1)$ is thus completely transparent, and in fact many
phenomena occuring in Grassmann algebras of higher rank cannot be
observed on such a simple example.

Nevertheless, the following result, which will be used later on,
shows that $\wedge(1)$ is an a sense `large enough' to form a target for
an onto homomorphism from every subalgebra of a Grassmann algebra
with a nontrivial odd sector.

\begin{lemma}
Let $A$ be a graded unital subalgebra of a finite dimensional Grassmann
algebra $\wedge(q)$.
If the odd part of $A$ is
nontrivial, then there is a surjective graded algebra homomorphism
$h\colon A\to\wedge(1)$.
\label{notorious}
\end{lemma}

\begin{proof}
Define $m$ as the smallest cardinality of a subset
$\beta\subseteq \{1,2,\dots,q\}$ such that in the expansion
(\ref{exp}) of at least one odd element $a\in A$ one has $a_{\beta}\neq 0$.
It follows from our assumptions on $A$ that $m$ is positive and odd.
Fix such a $\beta$.

Denote by $\zeta$ a fixed odd generator of the Grassmann
algebra $\wedge(1)$ and define a grading-preserving
$k$-linear mapping $\wedge(q)\to\wedge(1)$
by the rule
\begin{equation}
\sum_{\alpha\in 2^{\{1,\dots,q\}}} a_{\alpha}\xi^{\alpha}
\mapsto a_{\emptyset}+ a_{\beta}\zeta
\label{map}
\end{equation}
It is clear that the image of $A$ under the above mapping is all of
$\wedge(1)$, and thus it suffices to prove that its restriction to
$A$, say $h$, is an algebra homomorphism.
(It is easy to see that in general the linear map (\ref{map}) is not
a homomorphism on all of $\wedge(q)$.)

Let $x,y\in A$ be arbitrary.
Represent $x=x_0+x_1\xi^{\beta}+x_2$, $y=y_0+y_1\xi^{\beta}+y_2$, where
$x_0,y_0\in k$, and
in the expansion (\ref{exp}) of $x_2$ and $y_2$ both the constant terms
and the terms of order $\beta$ vanish. One has
\begin{eqnarray}
h(x) &=& x_0+x_1\zeta, \nonumber \\
h(y) &=& y_0+y_1\zeta,
\end{eqnarray}
and consequently
\begin{equation}
h(x)h(y)=x_0y_0+(x_0y_1+x_1y_0)\zeta
\end{equation}
At the same time,
\begin{eqnarray}
xy &=& (x_0+x_1\xi^{\beta}+x_2)(y_0+y_1\xi^{\beta}+y_2) \nonumber \\
&=& x_0y_0+(x_0y_1+x_1y_0)\xi^{\beta}+ \nonumber \\
& & x_0y_2+y_0x_2+x_1y_1\cdot 0+
(x_1y_2+x_2y_1)\xi^{\beta}+x_2y_2
\end{eqnarray}
We claim that the image under $h$ of the terms in the last line of
the above formula is zero, which finishes the proof.
Firstly, since $x_0,y_0\in k$ and $h(x_2)=h(y_2)=0$, it follows
from the linearity of $h$ that
$h( x_0y_2+y_0x_2)=0$. Both $x_2$ and $y_2$ have no constant term,
which implies that the term of order $\beta$ in the expansion
(\ref{exp}) of $(x_1y_2+x_2y_1)\xi^{\beta}$ vanishes
(as well as the constant term of course) and,
as a consequence, $h((x_1y_2+x_2y_1)\xi^{\beta})=0$.
Finally, $x_2y_2$ has no constant term and it cannot have non-vanishing
term of order $\beta$ either, because otherwise either $x_2$ or
$y_2$ would have contained a monomial of the form
$a_{\gamma}\xi^{\gamma}$ with $a_{\gamma}\neq 0$,
$\vert\gamma\vert$ odd, and $\vert\gamma\vert<m=\vert\beta\vert$,
which is impossible by the choice of $\beta$.
\end{proof}

In conclusion, notice that the choice of a system of free odd generators
in a Grassmann algebra is by no means unique, as there is no canonical
way to select it, but this non-uniqueness does not affect the concepts or
results above.

A very detailed treatment of Grassmann algebras from the viewpoint of
supergeometry is to be found in the book \cite{Be2}.

\section{Basic notions of supergeometry}
\subsection{}
Let $X$ be a topological space. Recall that the topology of $X$
forms a category, ${\frak T}(X)$, with inclusion mappings as morphisms.
A {\it presheaf} of graded algebras on
$X$ is a contravariant functor,
$\cal S\colon U\mapsto {\cal S}(U)$, from ${\frak T}(X)$
to the category of all graded algebras and grading-preserving
homomorphisms, with the requirement that $\emptyset\mapsto \{0\}$.
The elements of ${\cal S}(U)$ are called {\it sections} of
$\cal S$ over $U$. If $V\subseteq U$, then the image of the inclusion
$V\hookrightarrow U$ under $\cal S$ is denoted by
$\rho^U_V$ and called the {\it restriction morphism.} A common shorthand
for $\rho^U_V(f)$ is $f\vert_V$.
A presheaf, $\cal S$, on $X$ is called a {\it sheaf} if it satisfies
the following two axioms. Suppose $\gamma$ is an open cover of an
open subset $U\subseteq X$.
\smallskip

\begin{enumerate}
\item If $f\in{\cal S}(U)$ is such that $\rho^U_V(f)=0$ for
all $V\in\gamma$, then $f=0$.
\item Suppose there is a collection $f_V\in{\cal S}(V)$, $V\in\gamma$,
such that whenever $V,W\in\gamma$, one has $\rho^V_{V\cap W}(f_V)=
\rho^W_{V\cap W}(f_W)$. Then there is an $f\in{\cal S}(U)$ with
$\rho^U_V(f)=f_V$ for all $V\in\gamma$.
\end{enumerate}
\smallskip

A {\it stalk} of a sheaf on $\cal S$ at a point $x\in X$ is
the graded algebra direct limit
${\cal S}_x=\varinjlim\{{\cal S}(U)\colon U\ni x \}$.
If $f\colon X\to Y$ is a continuous mapping between topological spaces
and $\cal S$ is a sheaf on $X$, then the {\it direct image sheaf,}
$f_\ast({\cal S})$, on $Y$ is defined through 
$f_\ast({\cal S})(U)={\cal S}(f^{-1}(U))$ for all open $U\subseteq Y$.

For the basics of sheaf theory, the reader may consult e.g.
\cite{Gd}.

A {\it locally ringed superspace} (else: {\it geometric superspace})
{\it over} $k$ is a pair ${\frak X}=(X,{\cal S})$, where
$X$ is a topological space and $\Cal S_X$ is a sheaf of
local graded commutative $k$-algebras.
The latter means, through a slight abuse of
language, that for every $V\neq\emptyset$, ${\cal S}(V)$ is a
local graded commutative algebra, and that
 every stalk $\Cal S_{X,x}$ is a local graded commutative
algebra. The unique maximal ideal of the
unital algebra $\Cal S_{X,x}$ will be denoted by
${\frak m}_{X,x}$, and the corresponding residue field
$\Cal S_{X,x}/{\frak m}_{X,x}$ by $k_X(x)$ or simply
$k(x)$. Sections $f\in\Cal S(U)$
of the structure sheaf of a geometric superspace
are called {\it superfunctions} on an open set $U$.

A {\it morphism of geometric superspaces},
$f\colon {\frak X}\to {\frak Y}$, $f=(f_0,f^\sharp)$, is formed by a
continuous map $f_0$ between the underlying
 topological spaces $X$ and $Y$,
and a local morphism of sheaves of unital $k$-algebras
$f^\sharp\colon\Cal S_Y\to f_{0,\ast}(\Cal S_X)$; locality means that
for every $x\in X$ the induced (in an obvious way)
homomorphism between stalks,
$f_x^\sharp\colon \Cal S_{Y,f_0(x)}\to \Cal S_{X,x}$, satisfies the condition
$f_x^\sharp({\frak m}_{Y,f_0(x)})\subseteq{\frak m}_{X,x}$.
 Here $f_{0,\ast}(\Cal S_X)$ is the direct image sheaf on $Y$;
in general it is not a sheaf of unital algebras.

Let ${\frak X}=(X,{\cal S})$ be a geometric superspace, and let
$U$ be a non-empty
open subset of $X$. Then $(U,{\cal S}\vert_U)$ is a geometric superspace,
which one may call an {\it open geometric subsuperspace of} $\frak X$.

For a section $\varphi\in\Cal S_X(U)$ over an open set
$U\subseteq X$ and for any $x\in U$
one can define the {\it value} of $\varphi$ at $x$, denoted by
$\varphi[x]$, as the image under the augmentation homomorphism from
$\Cal S_{X,x}$ to the residue field
$k(x)=\Cal O_{X,x}/\frak m_x$. 
The necessity to introduce the functor of points
(cf. the next Section) is explained by the fact that
 superfunctions --- and
therefore morphisms between superspaces ---
are not uniquely determined by the collection of their values.
(A similar phenomenon occurs in algebraic geometry with the
structure sheaves of schemes with nilpotent elements.)

Notice that every locally ringed space, ${\frak X}=(X,{\cal S})$,
becomes a locally ringed superspace if one puts the trivial grading
on the algebras of sections, ${\cal S}(U)$, making them coincide with
their even parts and setting the odd parts equal to $(0)$.
We will call such superspaces {\it purely even,} or else
{\it bosonic.}

Every locally ringed superspace, ${\frak X}=(X,{\cal S})$, has
a reflection in the category of locally ringed spaces and their morphisms,
which we will denote by ${\frak X}_{even}$.
It admits a very transparent description: the underlying topological
space of ${\frak X}_{even}$ is $X$, and for every open $U\subseteq X$
the algebra of sections is just ${\cal S}(U)_0$. It is easy to see
that  ${\frak X}_{even}$ is indeed a locally ringed space
(a purely even locally ringed superspace). The pair formed by the
identity mapping of $X$ and the embedding of the sheaf
${\cal S}_0$ into ${\cal S}$ forms a superspace morphism from
$\frak X$ to ${\frak X}_{even}$, which we will denote by $r_{even}$.
It has the following universal property: for every purely even
locally ringed superspace $\frak Y$ and every superspace morphism
$f\colon{\frak X}\to {\frak Y}$ there exists a unique morphism of
locally ringed spaces $\bar f\colon{\frak X}_{even}\to {\frak Y}$ such
that $f=\bar f\circ r_{even}$.

Every locally ringed superspace, $\frak X$, has also a purely even
coreflection, denoted by ${\frak X}_{red}$ and called the
{\it reduced subsuperspace} of $\frak X$. There is a canonical
morphism $i\colon {\frak X}_{red}\to {\frak X}$ such that every
morphism from a purely even geometric superspace, $\frak Y$, to $\frak X$,
factors through $i$. The structure sheaf of ${\frak X}_{red}$ is
the quotient sheaf of ${\cal S}_X$ by the sheaf of ideals generated
by the odd sector of ${\cal S}_X$.

\begin{example} \rm A definition of graded (super) topological space was
given in
\cite{HM3}, and used to prove a `superised' Haar theorem (i.e., it was
proved that on a
graded topological group the only  Berezinian measure invariant
by graded translations is the one associated with the ordinary Haar measure of
the group). In this construction the structure sheaf is locally isomorphic
to the tensor
product of the structure sheaf of the underlying ordinary topological space
times an
exterior algebra. However, this concept is rather restrictive, and we will
not be examining it in what follows.
\qed\end{example}
\subsection{Supermanifolds} An important example of a locally ringed superspace
over $k=\R$ is provided by
a {\it graded domain} $U^{m, n}$ {\it of dimension} $(m,n)$,
where $m,n$ are natural numbers. Its underlying topological
space is an open domain, $U$, in an $m$-dimensional Euclidean space,
while the structure sheaf is isomorphic to the sheaf of germs
of infinitely smooth mappings from $U$ to the Grassmann algebra of
rank $n$. In other words, for every open subset $V\subset U$ the graded
algebra of superfunctions on $V$, ${\cal S}(V)$, is isomorphic to the
graded tensor product ${\cal C}^\infty(V)\otimes \wedge(n)\cong {\cal
C}^\infty(V,\wedge(n))$,
where ${\cal C}^\infty$ is the sheaf of smooth real-valued functions on
$U\subset\R^m$.

By definition one has an epimorphism of sheaves of $\R$-algebras
$\pi\colon{\cal S}\to{\cal C}^\infty$. One easily checks that for every
$x\in U$ the
maximal ideal of ${\cal S}_x$ is the inverse image of the maximal ideal of
${\cal C}^\infty_x$ under $\pi$. The kernel of $\pi$ coincides with the
sheaf of nilpotents of $\cal S$, which we denote $\cal N$. The quotient ${\cal
N}/{\cal N}^2$ turns out to be a free sheaf of rank $n$ on $U$, i.e., it is
isomorphic to $[{\cal C}^\infty]^n$.

If $n=0$, the definition of a graded domain is identical with that of a smooth
domain of dimension $m$.
At the other end there is the case $m=0$,
leading to a `purely odd' superspace  which we
will denote by $pt_n$. Its underlying topological space
is a singleton, $\{\ast\}$, while the constant structure sheaf
has $\wedge(n)$ as the algebra of global sections.
The $(0,0)$-dimensional superdomain $pt_0$ is just a singleton considered
as a trivial smooth manifold. Notice also that the purely odd
superspace $pt_q$ makes sense for an arbitrary field $k$ and not just
for $k=\R$. One can introduce the concept of the spectrum of an
arbitrary graded-commutative algebra $\Lambda$ very much in the same
fashion as it is being done in algebraic geometry for commutative
algebras, 
then the superspace $pt_q$ is exactly $\operatorname{Spec}\,\wedge(q)$
\cite{L2}.

If $U=\R^m$ and $n$ is a fixed natural number, the corresponding graded domain
is called a {\it Euclidean superspace} and is denoted $\R^{m,n}$.

A (real)
{\it smooth finite dimensional supermanifold (graded manifold)},
$\frak X$, of
dimension $(m,n)$ is a geometric superspace over the ground field
$k=\R$ that is locally isomorphic to an $(m,n)$-dimensional
graded domain. The underlying topological space of $\frak X$ is a
smooth manifold $X$ of dimension $m$.
Every superdomain is a supermanifold. Every smooth manifold is
at the same time a
supermanifold of purely even dimension of the form $(m,0)$.

Also in this case one has an epimorphism of sheaves of $\R$-algebras
$\pi\colon{\cal S}\to{\cal C}^\infty_X$, whose kernel is the nilpotent subsheaf
$\cal N$ of $\cal S$, and the quotient
${\cal E}={\cal N}/{\cal N}^2$ is a locally free
${\cal C}^\infty_X$-module of rank $n$, i.e., it is the sheaf of sections of
rank $n$ real, smooth vector bundle $E$ on $X$. Moreover (as it follows
from the
local isomorphism of $(X,{\cal S})$ with a  graded domain)
$\cal S$ is {\it locally} isomorphic to the exterior algebra sheaf of
${\cal E}$ (the sheaf of sections of the vector bundle $\wedge E$, whose
fibre at
a point
$x\in X$ is the exterior algebra $\wedge E_x$).

{\it Morphisms} between supermanifolds are just the
geometric superspace morphisms described above.
Thus, supermanifolds form a full subcategory of that of locally ringed
superspaces and their morphisms.

The category of supermanifolds has direct products \cite{Ko,HM1,BBH1}.

\subsection{Global structure of supermanifolds}
If $E$ is a rank $n$ vector bundle on an $m$-dimensional differentiable
manifold
$X$, and $\cal E$ is the sheaf of sections of $E$, let ${\cal S}=\wedge\cal E$
be the exterior algebra sheaf of ${\cal E}$, i.e.~the sheaf of sections
of $\wedge E$. It is quite easy to check
that $(X,{\cal S})$ is an $(m,n)$ dimensional supermanifold.
The vector bundle that, according to our previous discussion, can be associated
to  $(X,{\cal S})$, is  straightforwardly proved to be isomorphic to
$E$.

We may wonder whether this construction is general, in the sense that, given
a supermanifold $(X,{\cal S})$, the sheaf ${\cal S}$ is {\it globally}
isomorphic
to the exterior algebra sheaf of ${{\cal N}}/{{\cal N}}^2$. This is indeed
true, and
this is  usually known as {\it Batchelor's theorem}.
\begin{thm}\label{bat} Let $(X,{\cal S})$
be a supermanifold. The sheaves
${\cal S}$ and $\wedge({\cal N}/{\cal N}^2)$ are globally isomorphic as sheaves
of graded commutative $\R$-algebras.
\qed\end{thm}
This isomorphism is not canonical; as a matter of fact, the isomorphisms
between
${\cal S}$ and $\wedge({\cal N}/{\cal N}^2)$ are in a one-to-one
correspondences
with the sections of the epimorphism $\pi\colon{\cal S}\to{\cal C}^\infty_X$,
namely, with the morphisms of sheaves of $\R$-algebras $\sigma\colon{\cal
C}^\infty_X\to\cal S$ such that $\pi\circ\sigma=\o{id}_{{\cal
C}^\infty_X}$. The
original proof of Theorem
\ref{bat} involves nonabelian sheaf cohomology \cite{Ba0}.  A
deformation-theoretic proof was given by Blattner and Rawnsley \cite{BR}; a
detailed account of the latter  is given in
\cite{BBH1}.

The validity of Batchelor's theorem relies on the fact that the structure sheaf
$\cal S$ of a (real) supermanifold has trivial \u Cech
cohomology since
it admits partitions of unity. It is for instance known that Batchelor's
theorem does not hold for complex (holomorphic) graded manifolds,
cf.~\cite{Gre}.

Meticulous introductions to locally ringed superspaces and supermanifolds
(graded manifolds) are to be found in \cite{Be2,Ko,L1,L2,Ko,Ma1,BBH1}.

\section{Functor of points}
The functor of points traces its origins to algebraic geometry. Here we will
present it in the form it assumed in supergeometry.

Let $\frak X$ be an arbitrary superspace, and let $q\in\N$.
A $q$-{\it point} of $\frak X$ is any superspace
morphism $\varkappa\colon pt_q\to\frak X$.

We will first establish the following analogue of a
well-known result holding for locally ringed spaces \cite{DG}.

\begin{prop}
\label{o} Let ${\frak X}=(X,{\cal S})$
be a locally ringed superspace over an arbitrary field $k$.
The $0$-points
of $\frak X$ are in a natural one-to-one correspondence with those points
$x\in X$ having $k$ as their residue field (`smooth points').
\end{prop}

\begin{proof} Since every unital algebra homomorphism between fields
is an isomorphism, the image of $\{\ast\}$ under a $0$-point must have
$k$ as its residue field. On the other hand,
a morphism
$pt_0\to \frak X$ is uniquely determined by the choice of the underlying
mapping
$\{\ast\}\to X$, which is in turn given by selecting a point in $X$.
\end{proof}

The following observation helps to clarify the origin of the terminology.

\begin{corol} If ${\frak X}=(X,{\cal S})$ is a supermanifold, then
$0$-points of
$\frak X$ are in a natural one-to-one correspondence with the points of the
underlying smooth manifold of $X$.
\end{corol}

\begin{proof} The result follows from the isomorphism
${\cal S}_x\simeq{\cal C}^\infty_x\otimes\wedge(n)$ holding
for every $x\in X$.
Here ${\cal C}^\infty$ is the sheaf of $C^\infty$ functions on $X$, and $n$ is
the odd dimension of $\frak X$.
\end{proof}

Denote the collection of all $q$-points of $X$ by $pt_q(X)$.
The following is obvious from this definition.

\begin{prop}
\label{sub}
 Let $U$ be an open subsuperspace of a locally ringed
superspace $X$. Then for every $q\in\N$ the set $pt_q(U)$ forms
a subset of $pt_q(X)$ in a natural way. \qed
\end{prop}

\begin{remark} A superspace need not have $q$-points at all.
Using Proposition \ref{o}, it is easy to
construct a nontrivial geometric superspace
$X$ such that for every $q$, the set $pt_q(X)$ is empty,
see e.g.~a similar example in \cite{P0}.
\label{nopoints}
\qed\end{remark}

Let $\Cal G$ denote the category of all
finite dimensional Grassmann
algebras and grading-preserving unital algebra homomorphisms.
The opposite category, ${\cal G}^{op}$, is equivalent to the
category of all supermanifolds of the form
$pt_q$, where $q$ varies over $\N$.

Let $\varphi\colon \wedge(q)\to\wedge(p)$ be a morphism of graded
algebras. It determines a superspace morphism
$\varphi^\sharp\colon pt_p\to pt_q$ going in the opposite direction.
If now $\varkappa\colon pt_q\to{\frak X}$ is a $q$-point of a locally
ringed superspace $\frak X$, then the composition
$\varkappa\circ \varphi^\sharp$ is
a $p$-point of $\frak X$. Thus, $\varphi$ determines a mapping
\begin{equation}
\varphi({\frak X})\colon pt_q({\frak X})\to pt_p({\frak X})
\label{fun}
\end{equation}
having the form
\begin{equation}
 pt_q({\frak X})\ni\varkappa\mapsto \varkappa\circ \varphi^\sharp
\in pt_p({\frak X}).
\end{equation}
Using this observation, it is easy to verify that the
correspondence
\begin{equation}
\wedge(q)\mapsto pt_q(X),~~\varphi\mapsto\varphi({\frak X})
\end{equation}
from the category $\Cal G$ to the category
 $\Cal{S}ets$ of sets and mappings is a covariant functor.
It is of course representable by its very definition, with
$X$ as the representing object:
\begin{equation}
pt_q(X)=\operatorname{Hom}\,(pt_q,X)
\end{equation}

\begin{defin}
Denote by $\Cal{S}ets^{\Cal G}$ the category formed by
all covariant functors ${\frak X}\colon
\Cal G\to\Cal{S}ets$ and naturally
defined functorial morphisms between them.
Objects of this category, $\frak X$, are
called {\it virtual superspaces}.
To maintain consistency of our notation, we will denote the
image of $\wedge(q)$ under a functor $\frak X$ by
$pt_q({\frak X})$.
A morphism from a virtual superspace $\frak X$ to a virtual superspace
$\frak Y$ (a functorial morphism) is a collection of mappings
$f_n\colon pt_n(\X)\to pt_n(\Y)$, $n\in\N$,
commuting with mappings between the sets of points induced by
each morphism between Grassmann algebras: if
$h\colon\wedge(n)\to\wedge(m)$ is such a morphism, then it induces
mappings $h(\X)\colon pt_n(\X)\to pt_m(\X)$ and
$h(\Y)\colon pt_n(\Y)\to pt_m(\Y)$, and the requirement making
$f$ into a functorial morphism is that
\begin{equation}
f_m\circ h(\X)=h(\Y)\circ f_n
\end{equation}
\qed
\end{defin}

By assigning to every locally ringed superspace $\frak X$ the
virtual superspace of the form
 $[\wedge(q)\mapsto pt_q({\frak X})]$,
one obtains a functor from the category of locally ringed superspaces
and superspace morphisms to the category of virtual superspaces and
their morphisms.
Indeed, every superspace morphism $f\colon {\frak X}\to {\frak Y}$
gives rise to a collection of mappings
$f_q\colon pt_q({\frak X})\to pt_q({\frak Y})$ in a consistent way.
Here
\begin{equation}
f_q(\varkappa)=f\circ\varkappa\in pt_q({\frak Y})
\label{efkju}
\end{equation}
for every $q$-point $\varkappa$ of $\frak X$.

For more on the relationship between smooth supermanifolds and
the associated virtual superspaces, see e.g. \cite{BL2, V}.

\begin{example} The set $pt_q(\R^{m,n})$
of $q$-points of the $(m,n)$-dimensional
Euclidean superspace $\R^{m,n}$ is
 the even sector of the graded vector space
$\wedge(q)\otimes \R^{m\vert n}$, where $\R^{m\vert n}$ stands for the
graded vector space $\R^m\oplus\R^n$. To put it otherwise,
$pt_q(\R^{m,n})$ is the set of elements of the vector space
\begin{equation}
[\R^m\otimes\wedge(q)_0]\oplus [\R^n\otimes\wedge(q)_1],
\end{equation}
where of course $\wedge(q)_i$, $i=0,1$ denote the even and odd
sector of $\wedge(q)$, respectively. \qed
\label{point}
\end{example}

\begin{remark}
The image, $pt_\Lambda(\R^{m,n})$,
under the functor of points determined by $\R^{m,n}$ of one or
another `grassmannian' algebra $\Lambda$ forms a graded $\Lambda$-module,
$\Lambda^{m,n}$, which was routinely
accepted as the basic object of superanalysis by many theoretical and
mathematical physicists. The resulting approach to supergeometry is
known as the {\it DeWitt--Rogers approach,} cf. 
\cite{DW, R1, R2, JP}. The approach we are following here is known as the
{\it Berezin--Leites--Kostant approach,} cf. 
\cite{BL1, Ko, L1, L2, Be2, Ma3}.
A functorial link between the two
approaches to supermanifolds was pointed out by
Le\u\i tes \cite{L2} and (independently) A. Schwarz \cite{Sch1,Sch2},
and remains largely unexplored to date. A brief discussion can be found in
\cite{Pen1} and \cite{BBHP3}. See also a paper by Schmitt
\cite{Schm}, containing an excellent account of the functor of points
in supergeometry. An early reference is a Stokholm preprint by Bernstein
and Leites \cite{BL2}. Some nontraditional aspects of the functor of
points in infinite dimensional geometry are discussed in \cite{Pe3}.
Notice that, if $\Lambda$ is an infinite-dimensional `grassmannian'
algebra, then it usually carries a natural locally convex algebra
topology which has to be taken into account in the definitions;
the emerging subtleties may be disruptive for the expected pattern of
results, cf. \cite{BP}.
\qed
\end{remark}

\begin{remark}
Some virtual superspaces are represented by
actual geometric superspaces, while others are not.
Rather than constructing relevant examples now, we will wait till
a host of such virtual superspaces will appear in a natural way in later
parts of this article, cf. Remark \ref{nonrep}. \qed
\end{remark}

\begin{remark}
The category $\Cal{S}ets^{\Cal G}$, being a category of functors
to $\Cal{S}ets$, has a
certain additional structure making it a {\it topos}
in the sense of \cite{J},
that is, a nonstandard model of set theory, and in this role
it has already received some attention \cite{Y}.
Being a topos leads to the existence of a {\it transfer principle:}
every statement, $\phi$, about sets made in a certain language
can be translated into a statement, $\phi^\uparrow$,
 about virtual superspaces, and $\phi$ is true if and only if
$\phi^\uparrow$ is true.
Here supergeometry comes close to
topos theory, though no serious investigation
of the extent to what the classical results can be `automatically superised'
through the topos $\Cal{S}ets^{\Cal G}$ has been done so far.
In particular, the structure of the topos of virtual superspaces must
be investigated in much
great detail, and the first step is to understand
the expressive power of the language associated to the
topos of virtual superspaces.
\qed\end{remark}

Virtual superspaces can be considered as `shadows' of actual objects of
supergeometry, or sometimes as `blueprints' for those objects still to be
constructed. They are of little use in themselves.
The authors of \cite{BMFS} have
stressed that the progress in some areas of mathematical physics is
hampered by the fact that though some objects (say, supermoduli spaces)
admit a pretty clear interpretation through the functor of points,
there are known no `genuinely geometric' objects of supergeometry
representing them --- while such objects, and  not their `shadows,'
are exactly what one needs for work.

\section{Virtual topological superspaces}
We begin with an auxiliary construction.
Let ${\frak X}=(X,{\cal S})$ be a locally ringed superspace over
a field $k$, let $q\in\N$, and let $f$ be a superfunction on
$X$. For an arbitrary $\varkappa\in pt_q(X)$, the sheaf
morphism $\varkappa^{\sharp}$ is in essence a graded algebra
homomorphism from the stalk ${\cal S}_{X,\varkappa_0(\ast)}$ to
the Grassmann algebra $\wedge(q)$. (Here $\ast$ is the only element
of the topological space underlying $pt_q$.)
Denote by $\tilde f_{\varkappa}$ the germ of $f$ at the point
$\varkappa_0(\ast)$.

Let
\begin{equation}
f_q[\varkappa]=\varkappa^{\sharp}\left(\tilde f_{\varkappa}\right)
\end{equation}
This is an element of the algebra of global sections of $pt_q$,
which is isomorphic to the Grassmann algebra $\wedge(q)$.
As $\varkappa$ runs over the set of all $q$-points of $\frak X$, we thus
obtain a mapping
\begin{equation}
f_q\colon pt_q({\frak X})\to \wedge(q)
\end{equation}
Notice that
for $q=0$, the element $f_0[\varkappa]$ coincides with the value
of $f$ at the point $\varkappa_0(\ast)$, that is, the image of the germ
of $f$ under the augmentation homomorphism
${\cal S}_{X,\varkappa_0(\ast)}\to k$.

\begin{defin}
Let $k$ be a {\it topological} field, and let
${\frak X}=(X,{\cal S})$ be a locally ringed superspace over $k$.
For every natural number $q$, we define the {\it natural topology}
on the set $pt_q({\frak X})$ as the coarsest topology with the
following property: for every open subset $U\subseteq X$ and every
superfunction $f\in{\cal S}(U)$, the mapping
$f_q\colon pt_q({U})\to\wedge(q)$ is continuous with respect
to the subspace topology on $pt_q(U)$ and the
standard product topology the Grassmann algebra supports as a
finite dimensional vector space over $k$.
\end{defin}

Here is a convenient reformulation of the above definition.

\begin{prop} The space $pt_q({\frak X})$ equipped with the natural
topology is canonically homeomorphic to the direct limit
topological
space $\varinjlim_UU$, where $U$ runs over all open subsets of
$X$ ordered by natural inclusion and each $U$
is equipped with the coarsest topology making every
mapping $f_q\colon U\to \wedge(q)$ continuous,
$f\in {\cal S}(U)$.
\qed
\end{prop}

In particular, the above reformulation shows that the natural topology
is well-defined on sets of $q$-points
for every locally ringed superspace over an
arbitrary topological field.

\begin{example} If the structure sheaf on a `purely even' superspace
${\frak X}=(X,{\cal S})$
is a subsheaf of that of germs of continuous $k$-valued functions on
$X$, then the natural topology on $pt_0({\frak X})$ is
contained in that induced from $X$. In particular, if $\frak X$ is either
a Tychonoff topological space with the sheaf of germs of continuous
functions,
or a smooth manifold with the sheaf of germs of smooth real-valued
functions, then
$X$ coincides with the set of all $0$-points and the natural
topology on $X$ is identical with the initial topology. \qed
\end{example}

\begin{example}
\label{toppo}
 The natural topology on the set of $q$-points
$pt_q(pt_p)\cong (\wedge(q)_1)^p$ is easily shown to coincide with
the product topology. \qed
\end{example}

\begin{example} Let $X$ be an arbitrary topological space. We make it
into a (purely even) locally ringed superspace by equipping $X$ with
the sheaf of germs of continuous real-valued functions, endowing all
algebras of sections with trivial (purely even) grading.
It is well known and easily proved that the stalks, ${\cal S}_x$, of a
locally ringed space thus
defined admit no non-trivial $\R$-valued derivations, and consequently
the only homomorphism ${\cal S}_x\to\wedge(q)$ is that of
augmentation, $f\mapsto f(x)$. Consequently, for every
$q\in\N$, the set of $q$-points of $X$ admits a canonical bijection
with $X$ itself, and the topology on $pt_q(X)$ is the completely
regular replica
of the topology of $X$. In particular, if $X$ is Tychonoff, then
$pt_q(X)$ is canonically homeomorphic to $X$ itself for each $q$.
\label{triv}
\qed
\end{example}

\begin{lemma} For every graded algebra morphism
$\varphi\colon\wedge(p)\to\wedge(q)$,
the corresponding mapping $\varphi({\frak X})\colon pt_p({\frak X})
\to pt_q({\frak X})$ is continuous with respect to the natural
topologies on both spaces.
\end{lemma}

\begin{proof} Let an open subset $U\subseteq X$ and
$f\in{\cal S}(U)$ be arbitrary. Since the functions
of the form $f_q\colon pt_q(U)\to\wedge(q)$ determine the topology
on $pt_q({\frak X})$, it is enough to verify that the `pull-back'
of $f_q$ on $pt_p(U)$ by means of the mapping
$\varphi({\frak X})\colon pt_p({\frak X})
\to pt_q({\frak X})$ is continuous. In other words, it suffices to check
the continuity of the mapping
\begin{equation}
f_q\vert_{pt_q(U)}\circ \varphi({\frak X})\colon pt_p(U)\to \wedge (q)
\end{equation}
To this end, it is enough to notice that
$f_q\circ \varphi({\frak X})$ is the composition of the continuous
mapping $f_p\vert_{pt_p(U)}$ with the graded algebra homomorphism
$\varphi$ which is also continuous as a linear mapping between
finite dimensional spaces.
\end{proof}

The following is an immediate corollary.

\begin{prop} The correspondence
\begin{equation}
pt_q\mapsto pt_q(X),~~\varphi\mapsto\varphi({\frak X})
\end{equation}
is a covariant functor
from the category $\Cal G$ to the category
 $\Cal{T}op$.
\label{assig}
\end{prop}

\begin{defin} A covariant
functor from $\Cal G$ to  $\Cal{T}op$ will be 
called a {\it virtual topological
superspace.} The category whose objects are the  virtual topological
superspaces, and whose morphisms are the corresponding functor
morphisms, will be denoted by ${\cal T}op^{\cal G}$. 
\end{defin}

\begin{lemma} For every morphism $f\colon {\frak X}\to {\frak Y}$
between two locally ringed susperpaces and for every $q\in\N$ the
mapping $f_q\colon pt_q({\frak X})\to pt_q({\frak Y})$
defined by formula \ref{efkju} is continuous
with respect to the natural topologies.
\end{lemma}

\begin{proof} For each open subset
$U\subseteq Y$ and every superfunction $g\in {\cal S}_Y(U)$, the pull-back
$g_q\circ f_q$ coincides with $(f^\sharp(g))_q$,
where $f^\sharp(g)$ is an element of the algebra of sections
$f_{0,\ast}({\cal S}_X)(U)$, canonically isomorphic to
${\cal S}_X(f_0^{-1}(U))$,
 and is therefore continuous on
$pt_q(f_0^{-1}(U))$.
\end{proof}

The following is an immediate corollary.

\begin{prop} The assignment of a virtual topological superspace
 to every locally ringed superspace described in Proposition
\ref{assig} is functorial (in a covariant way).
\qed
\end{prop}


The following structural result is very simple yet useful.

\begin{prop}
\label{fib}
Let $\X$
be a virtual topological superspace. Then for every
natural $q$, the space $pt_q(\X)$ forms a fibration over
$pt_0(\X)$ in a canonical manner.
If $\frak X$ is [determined by] a locally ringed superspace, then
the fibre over $x$
is homeomorphic to $\operatorname{Hom}({\cal S}_x,\wedge(q))$.
\end{prop}

\begin{proof}
The augmentation homomorphism
$\beta\colon\wedge(q)\to k\cong\wedge(0)$ determines a superspace morphism
$\beta_{\bullet}\colon pt_0(\X)\to pt_q(\X)$. The image of
$\beta_{\bullet}$ under the functor $\X$ is a continuous mapping
and therefore supplies the desired canonical
fibration $\beta_{\bullet}(\X)\colon pt_q(\X)\to pt_0(\X)$.
The inclusion $\wedge(0)\cong k\hookrightarrow\wedge(q)$,
$\lambda\mapsto\lambda\cdot 1$ is a homomorphism of unital graded
algebras and therefore determines a superspace morphism
$i\colon pt_q\to pt_0$; one thus obtains a continuous mapping
$i(\X)\colon pt_0(\X)\to pt_q(\X)$. The obvious property
$\beta\circ i = i$ implies that $\beta(\X)\circ i(\X) = i(\X)$,
that is, $\beta(\X)$ is a retraction of $pt_q(\X)$ onto a subspace
homeomorphic to $pt_0(\X)$, and in particular all fibres are nonempty.

Now assume that $\X$ is determined by a locally ringed superspace,
which we will for simplicity denote with the same symbol
$\X=(X,{\cal S})$.
According to Proposition \ref{o}, $0$-points of $\X$ correspond to
those points $x\in X$ having $k$ as their residue field. It follows that
if $\varkappa\colon pt_q\to\X$ is a $q$-point, then $\varkappa_0(\ast)=x\in X$
is a $0$-point of $X$, while $\varkappa^\sharp$ can be thought of as an
arbitrary
graded algebra homomorphism from the stalk ${\cal S}_x$ to
$\wedge(q)$. Notice that $\varkappa_0(\ast)$ is exactly
$\beta(\X)(\varkappa)$. Therefore,
the collection of all $q$-points $\varkappa$ with
$\varkappa_0(\ast)=x$ coincides with the fibre of the canonical fibration
$\beta_{\bullet}(\X)\colon pt_q(\X)\to pt_0(\X)$ over the point $x$.
Another way to describe this fibre is as the collection,
$\operatorname{Hom}({\cal S}_x,\wedge(q))$, of all graded algebra
homomorphisms from the stalk ${\cal S}_x$ to $\wedge(q)$.
The proof is thus finished.
\end{proof}

\section{Non-compactness of locally ringed topological superspaces}

It is natural to call a virtual topological superspace,
$\frak X$, {\it compact}
if for every $q$, the topological space of $q$-points,
$pt_q({\frak X})$, is compact. In other words, those
compact superspaces residing
within a particular fragment of supertopology that we are considering now
form objects of the
category ${\cal C}omp^{\cal G}$. However, here we will show that
the only occurences of such a phenomenon are in a sense trivial, and
thus, informally speaking,
the phenomenon of compactness along the odd directions is never
observed among virtual topological superspaces.

First of all, we need to define what does it mean that a virtual
topological superspace determined by a locally ringed superspace
is nontrivial in the odd sector.

\begin{defin}
Let $\frak X$ be a locally ringed superspace. We say that the
virtual topological superspace determined by $\frak X$ is
{\it trivial in the odd sector} if for every $q\in\N$ and
$\varkappa\in pt_q({\frak X})$ the graded subalgebra
$\varkappa^{\sharp}({\cal S}_x)$ of $\wedge(q)$
has trivial odd sector.
\end{defin}

Here is an equivalent reformulation of the same concept. Recall that
$i\colon {\frak X}_{red}\to {\frak X}$ is
the canonical morphism from the reduced subsuperspace
(even co-reflection) of $\frak X$.

\begin{prop} Let $\frak X$ be a locally ringed superspace. The
corresponding virtual topological superspace is trivial in the odd
sector if and only if for every $q\in\N$, the continuous mapping
$i_q\colon pt_q({\frak X}_{red})\to pt_q({\frak X})$ is a
homeomorphism. Equivalently, the functor associating a virtual
topological superspace to $\X$ factors through the
even coreflection $\X_{red}$.
\qed
\end{prop}

Put loosely, this is the case where the odd sector of a superspace,
$\X$, tells us nothing about the topology on $q$-points that is not
already encoded in the even subsuperspace $\X_{red}$.

\begin{lemma}
\label{tych}
Let $\X$ be a locally ringed superspace, and let $x\in pt_0(\X)$
and $q\in\N$.
The restriction of the natural topology
to the fibre,
 $\operatorname{Hom}({\cal S}_x,\wedge(q))$, of the canonical fibration
$pt_q(\X)\to pt_0(\X)$ over $x$ coincides with the
topology induced from the Tychonoff topology on the infinite product
$\wedge(q)^{{\cal S}_x}$ under the embedding
$\operatorname{Hom}({\cal S}_x,\wedge(q))\hookrightarrow
\wedge(q)^{{\cal S}_x}$.
\end{lemma}

\begin{proof} The natural topology on the fibre
 $\operatorname{Hom}({\cal S}_x,\wedge(q))$, formed by all $q$-points
$\varkappa$ with $\varkappa_0(\ast)=x$, is, by the definition, the coarsest
topology making every mapping of the form
\begin{equation}
\varkappa\mapsto f[\varkappa]
\end{equation}
continuous, where $f\in {\cal S}(U)$, and $U$ is an arbitrary open
neighbourhood of $x$. Let $h_\varkappa$ be a homomorphism from the
stalk ${\cal S}_x$ to $\wedge(q)$ associated to $\varkappa$, then
the natural topology is the coarsest one making every mapping of
the form $\varkappa\mapsto h_\varkappa(\tilde f)$ continuous, where
$\tilde f$ is the germ of a superfunction $f$ as above.
This is exactly the topology of simple convergence on elements of
the stalk ${\cal S}_x$, that is, the topology induced on
$\operatorname{Hom}({\cal S}_x,\wedge(q))$ from
$\wedge(q)^{{\cal S}_x}$, as required.
\end{proof}

The following result shows that among virtual topological superspaces
determined by locally ringed superspaces, every compact superspace
is trivial in the odd sector, that is, it comes from a locally ringed space
rather than superspace.

\begin{thm} Let $k$ be an infinite topological field, and let
${\frak X}=(X,{\cal S}_X)$ be a locally ringed
topological superspace over $k$.
Assume that the
topological space $pt_1({\frak X})$ (with the natural topology)
is compact. Then the virtual topological superspace determined by
$\frak X$ is trivial in the odd sector.
\label{trivial}
\end{thm}

\begin{proof} Assume that
the virtual topological superspace determined by
$\frak X$ is non-trivial in the odd sector, that is,
there is a $q\in\N_+$ and a $q$-point of $\frak X$, $\varkappa$,
such that $A=\varkappa^{\sharp}({\cal S}_x)$ has nontrivial odd sector
as a graded subalgebra of $\wedge(q)$.
According to Lemma \ref{notorious}, there is a surjective morphism of graded
algebras $j\colon A\to\wedge(1)$.
Denote by $\beta\colon A\to k$ the restriction
of the augmentation homomorphism $\wedge(q)\to k$ to $A$.
It is clear that for every even element $a_0\in A_0$
one must have $j(a_0)=\beta(a_0)$. From here it follows that for
every value of the parameter $\lambda\in k$ the linear mapping
$j_{\lambda}\colon A\to \wedge(1)$ determined by
$j_{\lambda}(a_0+a_1)=\beta(a_0)+\lambda j(a_1)$, $a_i\in A_i$,
$i=0,1$, is a graded algebra homomorphism. For two different values
$\lambda_1\neq\lambda_2$, the homomorphisms $j_{\lambda_1}$ and
$j_{\lambda_2}$ are distinct. Every composition of the form
$j_{\lambda}\circ \varkappa^{\sharp}$ is a graded algebra homomorphism
from ${\cal S}_x$ to $\wedge(1)$, and therefore determines
a $1$-point of $\frak X$, and, moreover, an element of the fibre
of the fibration $pt_1({\frak X})\to pt_0({\frak X})$ at the point
$x$ (Proposition \ref{fib}).
We will denote such a $1$-point by $x_{\lambda}$.
For different values of $\lambda$, the points $x_{\lambda}$ are
different. It follows from Lemma \ref{tych} that the set
of all points $\{x_{\lambda}\colon\lambda\in k\}$ equipped with
the natural topology is canonically homeomorphic to the basic field
$k$. Indeed, in the topology of pointwise convergence,
 a net $x_{\lambda_\nu}$
converges to a point $x_{\mu}$ if and only if for every $z\in\wedge(q)$
the net $j_{\lambda_\nu}(\varkappa^{\sharp}(z))$ converges to
$j_{\mu}(\varkappa^{\sharp}(z))$, which is easily shown to be equivalent
to the convergence $\lambda_\nu\to\mu$.
Moreover, the set $\{x_{\lambda}\colon\lambda\in k\}$ is
readily verified to form a one-dimensional affine subspace 
in $\wedge(q)^{{\cal S}_x}$ and therefore is closed
with respect to the (locally convex Hausdorff)
Tychonoff product topology on $\wedge(q)^{{\cal S}_x}$.
Since an infinite topological field is never compact \cite{Wi},
it means that the fibre of $pt_1({\frak X})$ over $x$ is non-compact.
But it is closed in $pt_1({\frak X})$, which is a contradiction.
\end{proof}

\begin{remark}
\label{nonre}
Observe that the category of virtual topological superspaces
possesses `compactifications.' Suppose a virtual
superspace, $\X$, is `Tychonoff' in the sense that for each
$q$, the space $pt_q(\X)$ is Tychonoff. Define for each $n\in\N$
\begin{equation}
pt_n(\beta\X)\overset{def}=\beta(pt_n(\X)),
\end{equation}
where $\beta$ denotes, as usual, the Stone--\u Cech compactification.
Every continuous mapping $f\colon pt_n(\X)\to pt_m(\X)$
extends to a unique continuous mapping
$\tilde f\colon \beta(pt_n(\X))\to \beta(pt_m(\X))$, which enables one
to turn the correspondence $\wedge(n)\to \beta(pt_n(\X))$ into a
covariant functor and indeed a virtual topological superspace, containing
$\X$ as a virtual topological subsuperspace in a natural fashion.

Since every supermanifold $\X$ is `super-Tychonoff' in the sense
that the natural topology on each set $pt_n\X$ is Tychonoff, it
admits a nice compactification in the category of virtual
topological superspaces.
The compactification procedure for virtual topological superspaces
certainly deserves further attention. However, we want to stress that it does
not provide an answer, or at least a complete answer, to
 the problem of compactifying supermanifolds --- simply because,
as we will see shortly, some of the most intriguing hypothetical objects
of supertopology, such as the purely odd projective superspace, do not
correspond to any virtual topological superspace. The setting of functor
of points is, thus, too restrictive.
\end{remark}

\begin{remark}
\label{nonrep}
The above construction enables one, nevertheless, to produce 
numerous examples of virtual
topological superspaces that do not come from locally ringed superspaces.
Such is, for instance, the above described `compactification' of [the virtual
topological superspace assigned to] any
supermanifold, $\X$, whose odd dimension $n\neq 0$. 

Let us consider the simplest case, that of $\X=pt_1$. 
Assume that $\beta(pt_1)$ is of the form $\wedge(q)\mapsto pt_q(\Y)$ 
for some locally ringed superspace $\Y$. Notice that then
\begin{equation}
pt_0(\Y)\cong\beta(\{\ast\})\cong\{\ast\},
\end{equation}
 and
\begin{equation}
pt_1(\Y)\cong\beta(pt_1(pt_1))\cong
\beta(\wedge(1)_1)\cong\beta\,\R.
\end{equation}
According to
Theorem \ref{trivial}, one can assume without loss in
generality that $\Y$ is purely even, that is, a locally ringed
space. In such a case, all $1$-points of $\Y$ are just
$0$-points, that is, $pt_1(\Y)\cong\{\ast\}$, a contradiction. \qed
\end{remark}

\section{\label{g-spaces} Virtual topological superspaces as topological
$G$-spaces}

In this Section we will show how the theory of virtual topological
superspaces can be fused into the setting of abstract topological
dynamics.

Denote by $\wedge(\infty)$ the Grassmann algebra of infinite rank,
that is, the associative unital graded algebra freely generated by a
countably infinite set $\{\xi_1,\xi_2,\dots,\xi_n,$ $\dots\}$ of
pairwise anticommuting elements. (Cf. e.g. \cite{R2}.)
 We will fix a family of generators
in what follows.
The algebra $\wedge(\infty)$ can
be represented as the direct limit (in fact, union) of the increasing
family of Grassmann algebras $\wedge(n)$ of finite rank.
We will always assume that $\wedge(n)$ sits canonically inside
$\wedge(\infty)$ by identifying the former algebra
with the subalgebra of the latter generated
by the first $n$ free odd generators, $\xi_1,\xi_2,\dots,\xi_n$.
This standard embedding will be denoted by
$i_n\colon\wedge(n)\hookrightarrow\wedge(\infty)$.
We will denote by $\pi_n$ the
canonical retractive homomorphism from $\wedge(\infty)$ to
$\wedge(n)$, sending the generators
$\xi_i$ to themselves for $i=1,2,\dots,n$
and to zero for $i>n$. We will also denote the canonical embedding
of $\wedge(n)$ into $\wedge(m)$, $n\leq m$, by $i_{n,m}$.
Notice that $i_{m,n}=\pi_m\circ i_n$.

Denote by $\operatorname{End}\wedge(\infty)$ the (unital) semigroup of all
graded algebra endomorphisms of $\wedge(\infty)$.
It can be identified with the set
$\wedge(\infty)_1^\omega$ of all infinite sequences of odd
elements of $\wedge(\infty)$, because such sequences are in a
natural one-to-one correspondence with endomorphisms of $\wedge(\infty)$.
(This identification is not canonical though, as it depends on the
choice of a family of odd generators for $\wedge(\infty)$.)
Moreover, the semigroup operation on such sequences can be easily
interpreted in terms of substitution of variables.
(The interested reader may consult \cite{Be2} to see how it is being
done for finite-dimensional Grassmann algebras in the context of
change of odd variables and Berezin integration; the extension of
the procedure to the algebra of infinite rank is straightforward and
indeed insightful.
However, we are not going to make any use of it in this article.)

Let the symbol $\E$ stand for the subsemigroup of
$\operatorname{End}\wedge(\infty)$ consisting of all endomorphisms
$f\colon \wedge(\infty)\to\wedge(\infty)$ with finite dimensional
range. Equivalently, an $f\in\operatorname{End}\wedge(\infty)$
is in $\E$ if and only if one has $f(\wedge(\infty))\in \wedge(n)$
for a suitable $n\in\N$, that is, $f=i_n\circ\pi_n\circ f$.

Observe also that every $\pi_n$ belongs to $\E$, and moreover if
$f\in\operatorname{End}\wedge(\infty)$, then both
$f\pi_n$ and $\pi_nf$ are in $\E$.
 One can identify $\E$ with the set of
all elements of $\wedge(\infty)_1^\omega$ all of whose coordinates
depend (as polynomials) on the same finite collection of odd variables.

Let $\frak X$ be a virtual superspace, that is, a covariant functor
from $\cal G$ to ${\cal S}ets$. Denote by
\begin{equation}pt_{\infty}(\X)=
\varinjlim\{pt_q({\frak X}), pt_q(i_{n,m})\}
\label{direct}
\end{equation}
the limit of the direct system
 $\{pt_q({\frak X}), pt_q(i_{n,m})\}$
of sets and mappings.
Recall that elements of such a direct limit are equivalence classes
of the disjoint union $\sqcup_{q\in\N}pt_q({\frak X})$ under the
equivalence relation defined as follows: an element
$x\in pt_q(\X)$ is equivalent to an element $y\in pt_n(\X)$
if and only if for some $m\geq n,q$ one has
$pt_q(i_{q,m})(x)=pt_n(i_{n,m})(y)$.
The notation we will use for two equivalent elements: $x\sim y$.
The equivalence class of an element $x\in pt_n(\X)$ is denoted either
by $[x]$ or else by $\hat\imath_n(x)$,
and thus one obtains canonical mappings
$\hat\imath_n\colon pt_n(\X)\to pt_\infty(\X)$, $n\in\N$.

Every projection $\pi_n\colon\wedge(\infty)\to\wedge(n)$ determines
a mapping $\hat\pi_n\colon pt_{\infty}(\X)\to pt_n(\X)$ as the direct
limit of mappings $\widehat{(\pi_n\circ i_j)}\colon pt_j(\X)
\to pt_n(\X)$, $j\to\infty$.
The following result is readily deduced from the functoriality
of $\X$.

\begin{lemma} For every $n\in\N$,
$\hat\pi_n\hat \imath_n=\operatorname{Id}_{pt_n(\X)}$. \qed
\label{lee}
\end{lemma}

 If $\X$ is a virtual topological superspace,
the direct system in (\ref{direct})
consists of topological spaces and continuous
mappings, and the set
$pt_\infty(\X)$ carries a natural topology, making it into the
the topological space direct limit: the topology is
by definition the finest one making each mapping $\hat\imath_n\colon
pt_n(\X)\to pt_\infty(\X)$ continuous.

\begin{example}
Let $\R^{1,1}$ denote, as before, the standard $(1,1)$-dimensional
smooth superdomain. It follows from Examples \ref{point} and \ref{toppo}
that for every $q$ the topological space
$pt_q(\R^{1,1})$ is canonically homeomorphic with
the (underlying topological space of) the Grassmann algebra
$\wedge(q)$. Consequently, the topological space $pt_{\infty}(\R^{1,1})$
is homeomorphic with the injective limit of topological spaces
$\wedge(q)$, $q\to\infty$
under the natural inclusions, and thereby with the
(underlying topological space of) the infinite
dimensional Grassmann algebra $\wedge(\infty)$, equipped with the
finest locally convex topology. (As a topological
vector space, $\wedge(\infty)$ is
isomorphic to $\R^\omega$, equipped with the finest locally
convex topology, which topology in its turn
is well known to coincide, for the countable number of direct summands,
with the box product topology.) \qed
\end{example}

\begin{example} If, as in Example \ref{triv}, $X$ is a topological
space made into a purely even locally ringed superspace in a natural
way, then $pt_{\infty}(X)$ is canonically homeomorphic to $X$ itself.
\qed
\end{example}

\begin{prop}
The semigroup $\E$ acts in a natural way on the set $pt_\infty(\X)$.
If $\X$ is a virtual topological superspace, the action of
$\E$ upon the topological space $pt_\infty(\X)$ is one by
continuous transformations.
\label{action}
\end{prop}

\begin{proof}
Let $g\in\E$.
For every $n\in\N$, $g_{n}\overset{def}=g i_n$
is a graded algebra homomorphism
from $\wedge(n)$ to
$\wedge(\infty)$, and for a suitable $j$, one has
$g_n=i_j\pi_jg_n$.
The graded homomorphism
\begin{equation}
g_{j,n}\overset{def}=
\pi_jg_n\equiv\pi_jgi_n\colon\wedge(n)
\to\wedge(j)
\end{equation}
determines a mapping from $pt_n({\X})$ to $pt_j(\X)$,
which we will denote by $\hat g_{j,n}$.
Set $\hat g_n=\hat\imath_j\hat g_{j,n}$;
it is easy to see that the definition of the mapping
$\hat g_n\colon pt_n({\X})\to pt_{\infty}(\X)$
is independent of the actual
choice of $j$ as long as $\wedge(j)$ contains the range of
$g$.

If
$x\in pt_n(\X)$, $y\in pt_m(\X)$, and $x\sim y$, then quite evidently
$\hat g_n(x)=\hat g_m(y)$. This means that
the system of mappings
$\hat g_{n}$, $n\in\N$, gives rise to a mapping of $pt_{\infty}(\X)$ to
itself, which we will denote by the same letter $g$ as the original
element of the semigroup $\E$. This mapping is called the {\it motion} of
$pt_{\infty}(\X)$ determined by $g$, or simply a $g${\it -motion.}
For every $\kappa\in\wedge(\infty)$, one has
\begin{equation}
g(\kappa)=[\hat g_n(x)],
\label{class}
\end{equation}
where $n\in\N$, $x\in\wedge(n)$, $[x]=\kappa$, and the square brackets
symbolize the equivalence class containing a
given element. As it is customary in dynamics, we will often write $g\kappa$
instead of $g(\kappa)$.

Notice that in the topological case all participating mappings
are continuous, including the motion mapping $g\colon pt_{\infty}(\X)\to
pt_{\infty}(\X)$.

It remains to verify that for
each $\varkappa\in pt_\infty(\X)$, one has $(gh)(\varkappa)=g(h(\varkappa))$.
To prove this, choose an element, $x\in pt_q(\X)$, of the equivalence
class $\varkappa$ for some $q\in\N$.
Let $m,n\in\N$ be arbitrary natural numbers
such that $q\leq m\leq n$ and $\operatorname{range}(h_q)\subseteq
\wedge(m)$,  $\operatorname{range}(g_m)\subseteq
\wedge(n)$. Then $\pi_m\circ h\circ i_q\colon\wedge(q)\to\wedge(m)$ and
$\pi_n\circ g\circ i_m\colon\wedge(m)\to\wedge(n)$ are graded
algebra homomorphisms, whose composition coincides with
$\pi_n\circ( g\circ h)\circ i_q\colon\wedge(q)\to\wedge(n)$.
Since $\X$ is a covariant functor, one must have
\begin{equation}
\hat g_{n,m}\hat h_{m,q}=\widehat{(gh)}_{n,q},
\end{equation}
and consequently
$\hat \imath_n\hat g_{n,m}\hat h_{m,q}=\hat \imath_n\widehat{(gh)}_{n,q}$,
that is, $\hat g_{m}\hat h_{m,q}=\widehat{(gh)}_{q}$.
Since $\hat\pi_m\hat \imath_m=\operatorname{Id}_m$ by Lemma
\ref{lee}, $\hat g_m\hat h_q=\hat g_m\hat\pi_m\hat \imath_m
\hat h_{m,q}=\widehat{(gh)}_{q}$.
In particular,  $\hat g_m(\hat h_q(x))=\widehat{(gh)}_{q}(x)$,
and according to (\ref{class}),
\begin{eqnarray}
g(h(\kappa))&=&g[\hat h_q(x)]\nonumber \\
&=&[\hat g_m(\hat h_q(x))] \nonumber \\
&=& [\widehat{(gh)}_{q}(x)]\nonumber \\
&=& (gh)(\kappa),
\end{eqnarray}
finishing the proof.
\end{proof}

\begin{thm}
The category of virtual superspaces is isomorphic to a full
subcategory of the category of $\E$-sets.
\label{isosets}
\end{thm}

\begin{proof} First of all, we wish to turn the assignment
of a $\E$-set to every virtual superspace, described in Proposition
\ref{action}, into a covariant functor between the corresponding
categories.
 Let $\X$ and $\frak Y$ be two virtual superspaces, and let
$f\colon\X\to\Y$ be a morphism between them, that is, a collection
of mappings $f_n\colon pt_n(\X)\to pt_n(\Y)$, $n\in\N$, satisfying
the requirement that whenever
$h\colon\wedge(n)\to\wedge(m)$ is a morphism of graded unital algebras,
then  $f_m\circ h(\X)=h(\Y)\circ f_n$, where
$h(\X)\colon pt_n(\X)\to pt_m(\X)$ and
$h(\Y)\colon pt_n(\Y)\to pt_m(\Y)$ are mappings images of $h$ under
the functor $\X$.

The rule
\begin{equation}
f(\kappa)=[\hat\imath_n f_n(x)]\label{isafunctor}
\end{equation}
determines a mapping from $pt_{\infty}(\X)$ to $pt_{\infty}(\Y)$,
and it follows directly
from the fact that $f$ is a functor morphism that our
newly-defined mapping commutes
with the action of the semigroup $\E$. The verification that the
assignment (\ref{isafunctor})
satisfies the functorial properties is easy, and thus our first objective is
achieved.

The next step is to show that different virtual superspaces are
being sent to different $\E$-spaces. To this end, note that
an arbitrary virtual supserspace, $\X$, can be fully recovered from
the $\E$-set $pt_{\infty}(\X)$. Indeed, it is enough to
observe that for every $n\in\N$, the set $pt_n(\X)$ is obtained as
the image of the mapping $\hat\varpi_n\colon pt_{\infty}(\X)\to
pt_{\infty}(\X)$, while for every morphism $h\colon\wedge(n)\to\wedge(m)$
the mapping $h(\X)\colon pt_n(\X)\to pt_m(\X)$ determined by it is
obtained by setting, for each $x\in pt_n(\X)$,
$h(\X)(x)=\hat\varpi_m(\widehat{i_m\circ h\circ\pi_n}(x))$.
Consequently, we obtain an isomorphism of the category of virtual
superspaces with a subcategory of $\E$-sets.

It remains to show that this subcategory is full, that is, we do not
get any additional morphisms on top of those
determined by morphisms between virtual superspaces.
With this purpose, observe that
every morphism of $\E$-sets, $f\colon pt_{\infty}(\X)\to
pt_{\infty}(\Y)$, determines a collection of mappings
$f_n\colon  pt_n(\X)\to pt_n(\Y)$ defined for each $n$ by an obvious
rule $f_n=   \hat\varpi_nf\hat\imath_n$.
Being equivariant, $f$ commutes with the mappings of the form
$\hat h_{m,n}$, where $h\colon\wedge(n)\to\wedge(m)$ is an arbitrary
morphism of Grassmann algebras, and this implies easily that the
collection $(f_n)_{n\in\N}$ is a functorial morphism. The morphism of
$\E$-sets determined by it is exactly $f$. This finishes the proof.
\end{proof}

The following result is obtained from the previous one word for word
by keeping track
of continuity of all participating mappings throughout the proof.

\begin{thm}
The category of virtual topological superspaces is isomorphic to
a full subcategory of the category of topological $\E$-spaces.
\label{isospaces}
\qed\end{thm}

\begin{remark} The subcategories that we obtain are proper in both
cases. For example, no $\E$-set, $X$, endowed with a constant action of
$\E$ of the form $gx=x_0$ for all $g\in\E$, $x\in X$, and a fixed
$x_0\in X$, is an image of a virtual superspace: every $\E$-set, $X$,
obtained from a virtual superspace is a union of its subsets of the
form $\hat\varpi_n(X)$.

If one wants to establish an isomorphism of the category of virtual
(topological)
superspaces with that of {\it all} (topological) $G$-spaces, then it
probably makes sense to consider the semigroup of {\it all}
continuous endomorphisms of $\wedge(\infty)$ (or other locally
convex Grassmannian
algebra, such as the DeWitt algebra $\Lambda_{\infty}$ with its
canonical Fr\'echet topology, see \cite{Pe2,BP} for examples),
so that the semigroup of transformations will have both
unity and topology, enabling one to
exercise some control over the behaviour of subsets of $X$
of the form $\hat\varpi_n(X)$. \qed
\end{remark}

\begin{remark} Notice that in both Theorems \ref{isosets}
and \ref{isospaces} we could only
speak of categories being isomorphic rather than equivalent:
indeed, an isomorphism between categories in question depends on the
choice of a set of free odd generators for the infinite-dimensional
Grassmann algebra $\wedge(\infty)$, and therefore
there is no canonical isomorphism in sight. \qed
\end{remark}

\begin{remark}
The procedure of compactification of virtual topological superspaces,
outlined in Remark \ref{nonre},
 can be easily
described in dynamical terms. Let $\X$ be an arbitrary virtual
topological superspace.
 Since the acting semigroup, $\E$, is
discrete, its action extends to the Stone--\u Cech compactification,
$\beta (pt_{\infty}(\X))$, of the topological $\E$-space
$pt_{\infty}(\X)$ simply by extending each motion
$g\colon pt_{\infty}(\X)\to pt_{\infty}(\X)$ to a continuous self-mapping
of the Stone--\u Cech compactification. The resulting $\E$-space,
$\beta (pt_{\infty}(\X))$, is easily seen to contain the $\E$-space
$pt_{\infty}(\beta\X)$ as an everywhere dense $\sigma$-compact
$\E$-subspace, where $\beta\X$ is the virtual compact space described
in Remark \ref{nonrep}. In general, $pt_{\infty}(\beta\X)$ is non-compact,
but rather a $k_\omega$-space (an easy example to check is
$\X=\R^{1,1}$). It means that the $\E$-space $\beta (pt_{\infty}(\X))$
does not, in general, come from a virtual topological space.
The $\E$-space $pt_{\infty}(\beta\X)$
is obtained as the union of closures formed in in $\beta (pt_{\infty}(\X))$
of all subspaces of the form $\hat\imath_n(pt_n(\X))$. \qed
\label{cco}
\end{remark}

The previous Remark \ref{cco} makes it evident that an
abstract dynamical approach, in
which the basic object of study is an arbitrary topological
$\E$-space, is even more general than the functorial one.
There is a certain room for theory development,
including formulating and proving analogues of all the major
classical results of topology.

However, it appears that even this approach to the problem is too
narrow to incorporate some of the much desired but as yet non-existent
objects of supertopology such as the purely odd projective superspace
(cf. the next Section.)

Nevertheless, it might well happen that the dynamical approach offers
a good vantage point from which to survey the present
state of the theory and map out future directions.

\section{Final discussion: compactness vs cohomology}
Now the reader is prepared to face the challenging concrete problems
dealing with the existence of nontrivial analogue of
topologies for superspaces, and in
particular, about  compactifications `along the odd directions'.
We believe that the best way to outline them is to quote directly
from three esteemed experts in the area, adding some
minimalistic comments of our own.

D.~A.~Leites was probably the first to put the problem forward. Here is
how he describes it in his problem survey article \cite{L3}, pp. 650--651.

\begin{quote}\small
(c) Everybody knows the importance of orbits of group actions, e.g.~those
in the
coadjoint representation host all the classical mechanics. Now the category of
supermanifolds is not closed with respect to supergroup actions. Consider for
instance $GL(n)$ acting on the space of the identity representation. There are
two orbits: the origin and the rest. If we now look at the space as an
$(0,n)$-dimensional supermanifold we see that the complement to the origin is
just a kind of halo, indescribable except in the language of the point functor.
The functor corresponding to the complement of the origin is not presented by a
supermanifold.

Functors on the category of commutative (super)algebras represented by
manifolds or supermanifolds are good because you can construct differential
or ``at least'' algebraic geometries on them.
{\it How to distinguish subfunctors corresponding to the orbits of
supergroup action?} (A similar problem takes place for groups and their
orbits in prime characteristics.)

{\it Is it possible to construct mechanics on such orbits, i.e. integral
and differential calculus?}
\end{quote}

A moment's thought shows that the
situation is even less favourable than it appears from
the above quotation:
 the hypothetical object of supergeometry described by
Leites as the principal orbit of the action of $GL(n)$ on the
purely odd dimensional supermanifold
$pt_n$ does not even correspond to a virtual superspace.
More exactly, the `functor corresponding to the complement of the
origin' referred to by Leites does not exist.
Assuming such a virtual superspace, $\X$, existed, one would obviously
have for $q=0$
\begin{equation}
pt_0(\X)=pt_0(pt_q)\setminus\{0\}=\emptyset
\end{equation}
Since $\X$ is a covariant functor, the above observation leads to that
\begin{equation}
pt_q(\X)\subseteq\beta(\X)^{-1}(pt_0(\X))=\emptyset
\end{equation}
for all $q$,
where $\beta\colon\wedge(q)\to\wedge(0)\cong k$ is the homomorphism
of augmentation.

This remark shows in a most striking way 
that none of the known frameworks for
supergeometry, not even the functor of points, allows for the existence of
some of the most interesting objects
one would like to see implemented in supertopology.

Here is how Yu.~I.~Manin has formulated the problem in his
survey paper \cite{Ma1}.

\begin{quote}\small
{\it How do we compactify a supermanifold in the odd directions?}

Apparently, the lack of this procedure hinders the construction of a
cohomology theory of supermanifolds in which the Schubert supercells
have classes that depend not only on their substructure.

Le\u\i tes put forward the conjecture that in a suitable category there
must be an object ``purely odd projective superspace'': the factor
of the complement
$\operatorname{Spec}\, k[\xi_1,\dots,\xi_n]\setminus\operatorname{Spec}
\,k$ under the action of the multiplicative group
$(t,\xi_i)\mapsto (t\xi_i)$. In the usual sense of the word, the
corresponding space is empty.
\end{quote}

In our notation, $\operatorname{Spec}\, k[\xi_1,\dots,\xi_n]=
pt_n$, while $\operatorname{Spec}\, k=pt_0$. Again, one cannot
associate to the
`purely odd projective superspace' a non-degenerate
virtual superspace for the same reason as before.

Ivan Penkov argues in \cite{Pen2} that a
satisfactory cohomology theory can hardly be constructed for general
supermanifolds without a good understanding of the compactification
in the odd
dimensions, and suggests a paradoxical idea that for certain
supergroups $G$
the right analogue of such a compactification procedure is
specifying an action of $G$ on an affine superspace. We quote \cite{Pen2}:
\begin{quote}\small
However, it is essential to note that, to our mind, the problem ... is
intimately bound with another two most important and interrelated problems in
the theory of complex (algebraic or analytic) supermanifolds; the problem of
compactification along odd directions, as formulated by Manin \cite{Ma1}, and
the problem of constructing a cohomology theory on supermanifols, in which the
cohomology groups are indexed by elements of $\Z_2[\varepsilon]$,
i.e., by elements of
the dimension semigroup of the supermanifolds. (Incidentally, only solutions to
these problems will yield a definite definition of the concepts of
``compactness'' and ``cohomology theory'' in question.)
\end{quote}
\centerline{......}
\begin{quote}\small
Thus, it is hardly likely that the required cohomology theory for general
supermanifolds could be constructed without a good understanding of odd
compactification. However, on supermanifolds of the form $G^0/P$ (which to all
appearances admit a canonical compactification along odd
directions),\footnote{Here $G^0$ is a classical Lie supergroup \cite{Pen1}, and
$P$ is a parabolic sub-supergroup [Note of the authors].} it should be
possible to
reach the goal directly.
\end{quote}

Discussing this
idea at any depth would lead us far astray from the topic of this paper,
and as far as we are aware, no further progress in this direction has
been made yet.

Instead, we conclude this paper with a few comments about cohomology theories
on supermanifolds. Instances of such theories were considered
since the very first developments of supergeometry. A cohomology theory
which is
quite natural to take into account is the de Rham cohomology of the graded
differential forms
\cite{Ko}. Let
${\frak X}=(X,{\cal S})$ be a supermanifold of dimension $(m,n)$, and let
${\D}er\,{\cal S}$ be the sheaf of graded $\R$-linear derivations of $\cal S$:
for every open set $U\subset X$, ${\D}er\,{\cal S}(U)$ is the space
(actually, a  graded ${\cal S}(U)$-module) of
all
$\R$-linear morphisms $D\colon{\cal S}(U)\to{\cal S}(U)$ satisfying the graded
Leibniz rule
$$D(fg)=D(f)\,g+(-1)^{\tilde D\,\tilde f}\,f\,D(g)$$
(whenever $D$ and $f$ are homogeneous; a tilde denotes here again the
grading of an element in a graded module). By a straightforward generalization
of what happens in  ordinary differential geometry, the sheaf
${\D}er\,{\cal S}$
plays here the role of the sheaf of sections of the tangent bundle to $\frak
X$. The graded $\cal S$-dual of ${\D}er\,{\cal S}$, denoted
$\Omega_{\frak X}^1$, is the sheaf of {graded differential 1-forms.} Now,
proceeding as in the ordinary case, one introduces for every natural number
$p$ the sheaf of graded differential $p$-forms $\Omega_{\frak X}^p$, and an
exterior differential $d\colon \Omega_{\frak X}^p\to\Omega_{\frak X}^{p+1}$
such that $d^2=0$. The {\it super de Rham cohomology of $\frak X$} is the
cohomology of the complex
$${\cal S}(X) \,\stackrel{d}{\longrightarrow}\,
\Omega_{\frak X}^{1}(X)\,\stackrel{d}{\longrightarrow}\,
\dots \,\stackrel{d}{\longrightarrow}\,
\Omega_{\frak X}^{p}(X)\,\stackrel{d}{\longrightarrow}\,
\Omega_{\frak X}^{p+1}(X)\,\stackrel{d}{\longrightarrow}\,\dots$$
(here $\Omega_{\frak X}^{p}(X)$ is the graded $\R$-vector space
of global graded $p$-form on $\frak X$, etc.)
If $\eta$, $\tau$ are {\it odd} (in the $\Z_2$-gradation) graded differential
1-forms, then $\eta\wedge\tau=\tau\wedge\eta$; as a consequence, the sheaves
$\Omega_{\frak X}^p$ are nontrivial for every $p$, i.e., the super de Rham
complex is not bounded from above. However this is the source of little
trouble,
in that {\it the super de Rham cohomology of $\frak X$ is
isomorphic to the ordinary de Rham cohomology of the differentiable manifold
$X$} \cite{Ko,BBH1}.

Recently a cohomology for supermanifolds has been proposed
\cite{V2} (related papers are also \cite{BS,BL6, BL4,VZ1,VZ2}),
which is claimed to be nontrivial, i.e.~to be in general
nonisomorphic to the de
Rham cohomology of the underlying manifold. However, it has also been claimed
that in order to be consistent this theory must be constrained in such a way
that it reduces once more to the super de Rham cohomology
theory above described
\cite{Vic1} (cf.~also \cite{Vic2} for related work, where still another
cohomology theory is proposed.)

Cohomology theories have also been considered in  a category of supermanifolds
which is wider than the one here considered. This category was introduced
in \cite{BB1}, and thoroughly studied in the papers \cite{BBHP1,BBHP2}
and in \cite{BBH1}, to overcome some inconsistencies of the original approach
by Rogers (\cite{R1,R2}; see also \cite{DW,Rt1}).
In particular a generalization of the ordinary de Rham
cohomology has been formulated for this category
\cite{Ra,BB2,Br2,BBH1}. Given a
supermanifold $(X,{\cal A})$ in this category, this cohomology is nonisomorphic
to the de Rham cohomology of the underlying differentiable manifold
essentially when
the sheaf $\cal A$ has nontrivial \v Cech cohomology (it should be noticed that
the \v Cech cohomology of the structure sheaf of the supermanifolds we have
considered in this paper is always trivial, as it happens with the ordinary
differentiable manifolds). Therefore this cohomology is sensitive to the
`superdifferentiable' structure of the supermanifolds rather than to its
`supertopology' (whatever this may mean), as it is shown by the examples
in \cite{BB2,BBH1}.

The study of this category of
supermanifolds falls beyond the scope of the present work.
Anyway, it remains unclear to what extent this category supplies
answers to some of the above problems, and to what extent merely seeks
to reformulate them in disguise.

\medskip
\noindent{\bf Acknowledgements.} We thank C.~Bartocci for valuable
discussions. The research the present paper is based upon has been
supported by SISSA, by the National Group for Mathematical Physics of C.N.R.,
by the Research Project `Geometria delle variet\`a differenziabili'
(all Italy), by the Internal Grants Committee of the Victoria
University of Wellington, and by
the Marsden Fund grant for basic research
 `Foundations of supergeometry' (both New Zealand).

\oskip\oskip\footnotesize
\noindent {\it 1991 Mathematics Subject Classification:}
Primary 54A99, 58A50, secondary 18F99, 54H20.
\par\noindent
{\it Keywords:} Locally ringed superspaces, virtual superspaces,
topological superspaces, topological spaces of points.

\end{document}